\documentclass[11pt]{amsart}   
\usepackage{graphics}                  

\newtheorem{theorem}{Theorem}

\theoremstyle{definition}

\theoremstyle{remark}

\numberwithin{equation}{section}

\newtheorem{conjecture}[theorem]{Conjecture}
\newtheorem{ga}[theorem]
{Genericity Assumption}

\def \~{\, \tilde {}\,}

\def \w{\widetilde}

\def \0{\emptyset}

\begin{document}


\title [On the number of inscribed squares in a simple closed curve ...]%
{On the number of inscribed squares in \\
a simple closed curve in the plane}




\author{Strashimir G. Popvassilev}

\address{The City College of New York, 160 Convent Avenue, New York, NY 10031}
\email{strash.pop@gmail.com}

\date{October 26, 2008}


\begin{abstract}
We show that for every positive 
integer $n$ there is a simple closed
curve in the plane (which can be 
taken infinitely differentiable 
and convex) which has exactly $n$ 
inscribed squares.  
\end{abstract}

\maketitle

\section*{Introduction}

It is an open problem if for every 
simple closed curve in the plane 
there are four points from the curve 
that form the vertices of a square. 
Such a square is called inscribed in 
the curve (though it is not required 
that it is contained in the region 
bounded by the curve). The problem is 
simply stated, old, and has only 
partial positive solutions. 
See \cite{Pak} for a list of papers, 
and for comments. 

The present note answers in the negative 
what we interpret as a conjecture 
posed by Jason Cantarella on his 
web page \cite{web}. The web site 
comments on his joint work with 
Elizabeth Denne and John McCleary 
on this problem. Their results have 
been announced in \cite{abs}. The 
author has recently been informed 
by Elizabeth Denne and Jason Cantarella 
that the preprint presenting the 
results announced in \cite{abs} 
is not yet ready to be released. 
Our recent discussion 
with Jason Cantarella and  
Elizabeth Denne on some of 
the ideas presented in \cite{web} 
and \cite{Denne}
has been helpful to the author, yet 
the following statement made at 
the web site \cite{web} has not been yet 
clarified: 

`Our results prove that there are an odd number of squares in any simple closed curve which is differentiable or ``not too rough''.'

Apparently the exact statement of the 
above result would appear in the forthcoming paper by Jason Cantarella, 
Elizabeth Denne and John McCleary. 

The purpose of the present note is the 
proof of the following. 

\begin{theorem} \label{exactly-n}
For every positive 
integer $n$ there is a simple closed
curve in the plane (which can be 
taken infinitely differentiable 
and convex) which has exactly $n$ 
inscribed squares.  
\end{theorem}

This seems to indicate (though we provide some ``evidence'' only, and no complete 
proof) that the following 
conjecture about the number of inscribed 
squares of an immersed in the plane curve (self intersections allowed) made at the same web site, is not valid, if only differentiability is assumed: 

`We might guess that the number of squares is equal to $St + (J^+  - J^- ) + 1$~mod~$2$.' 

As indicated in \cite{web}, 
$St$, $J^+$ and $J^-$ denote the invariants 
of the curve called strangeness, positive jump, and 
negative jump, introduced by Arnold. 
See \cite{Arnold}. 

\section{How to control the number of 
inscribed squares} 

First we sketch the construction of an 
infinitely differentiable simple closed 
curve in the plane that has exactly two 
inscribed squares.

Clearly the unit circle has infinitely 
many inscribed squares. On the other 
hand it is easy to modify the unit 
circle to obtain a (non-differentiable) 
simple closed curve which has exactly 
two inscribed squares. The construction is shown on Figure~1, left. 
The arc determined 
by central angles $5\pi\over4$ and $7\pi\over4$ is removed from the unit circle, and replaced by the semi-circle 
$y=\sqrt{{1\over2} - x^2} - 
{1\over\sqrt2}$. The reader may verify 
that there are only two inscribed squares, as shown on Figure 1, left.

\begin{figure*} 
 \scalebox{.16}{\includegraphics*{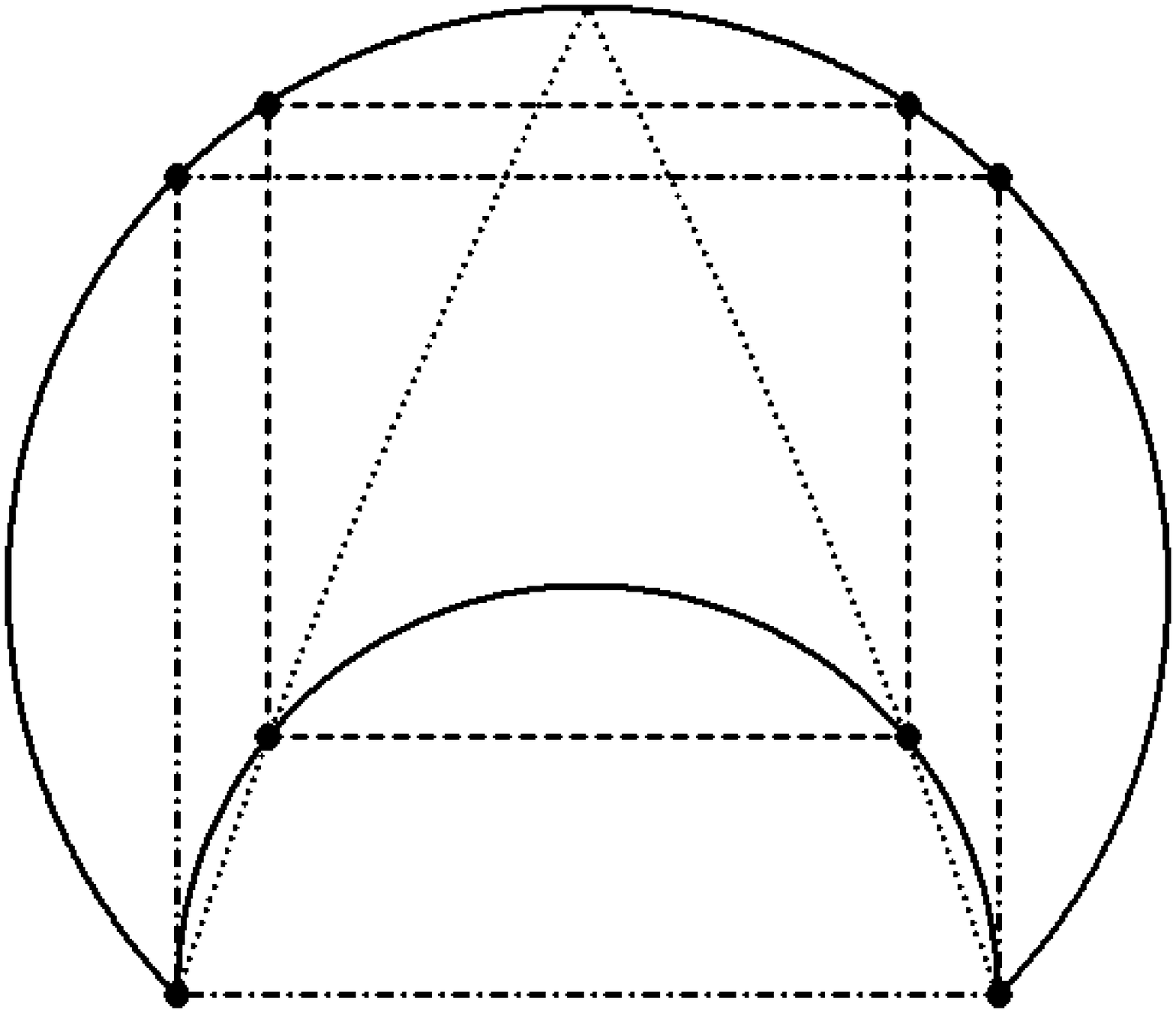}} \scalebox{.16}{\includegraphics*{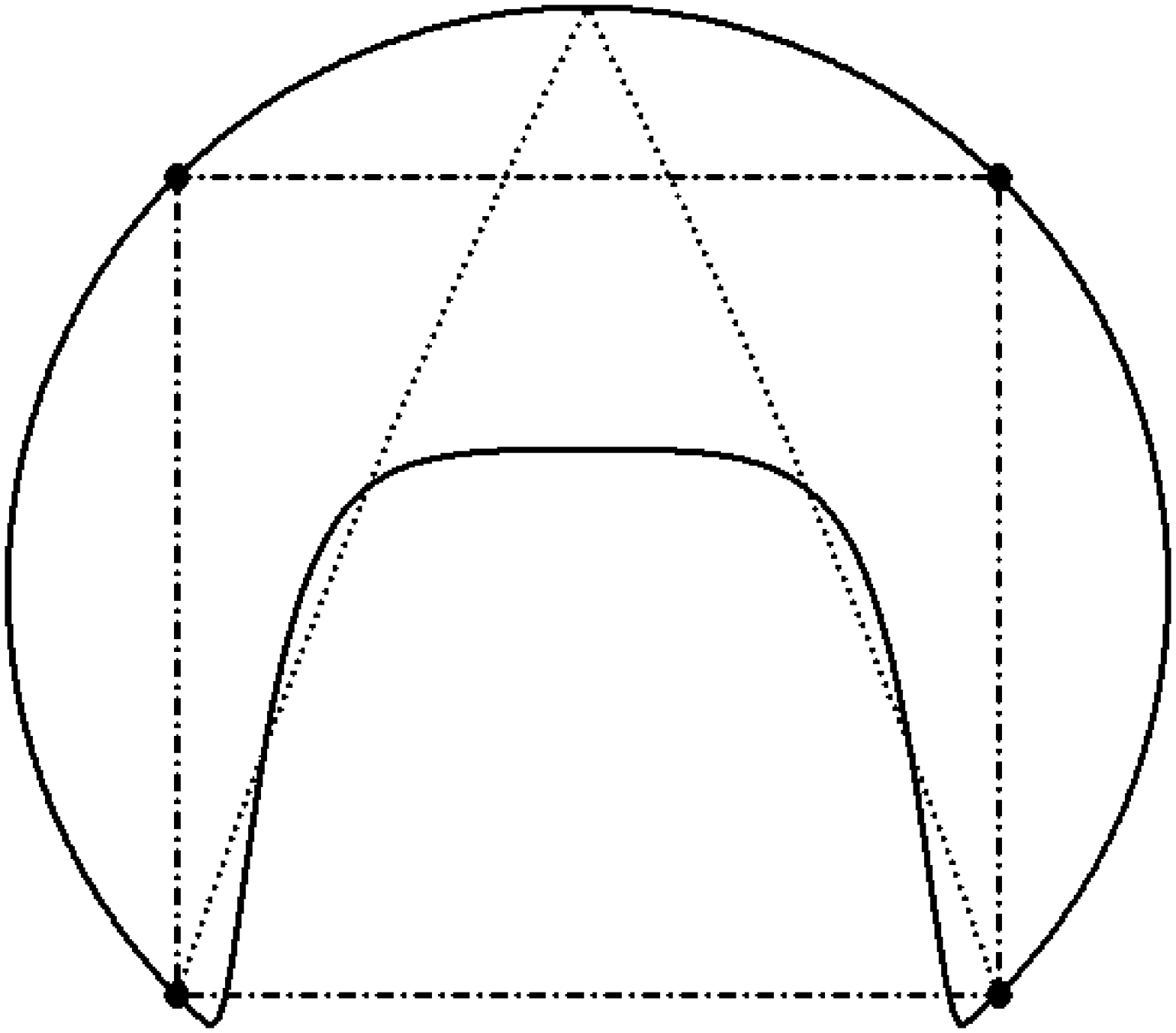}}
\scalebox{.16}{\includegraphics*{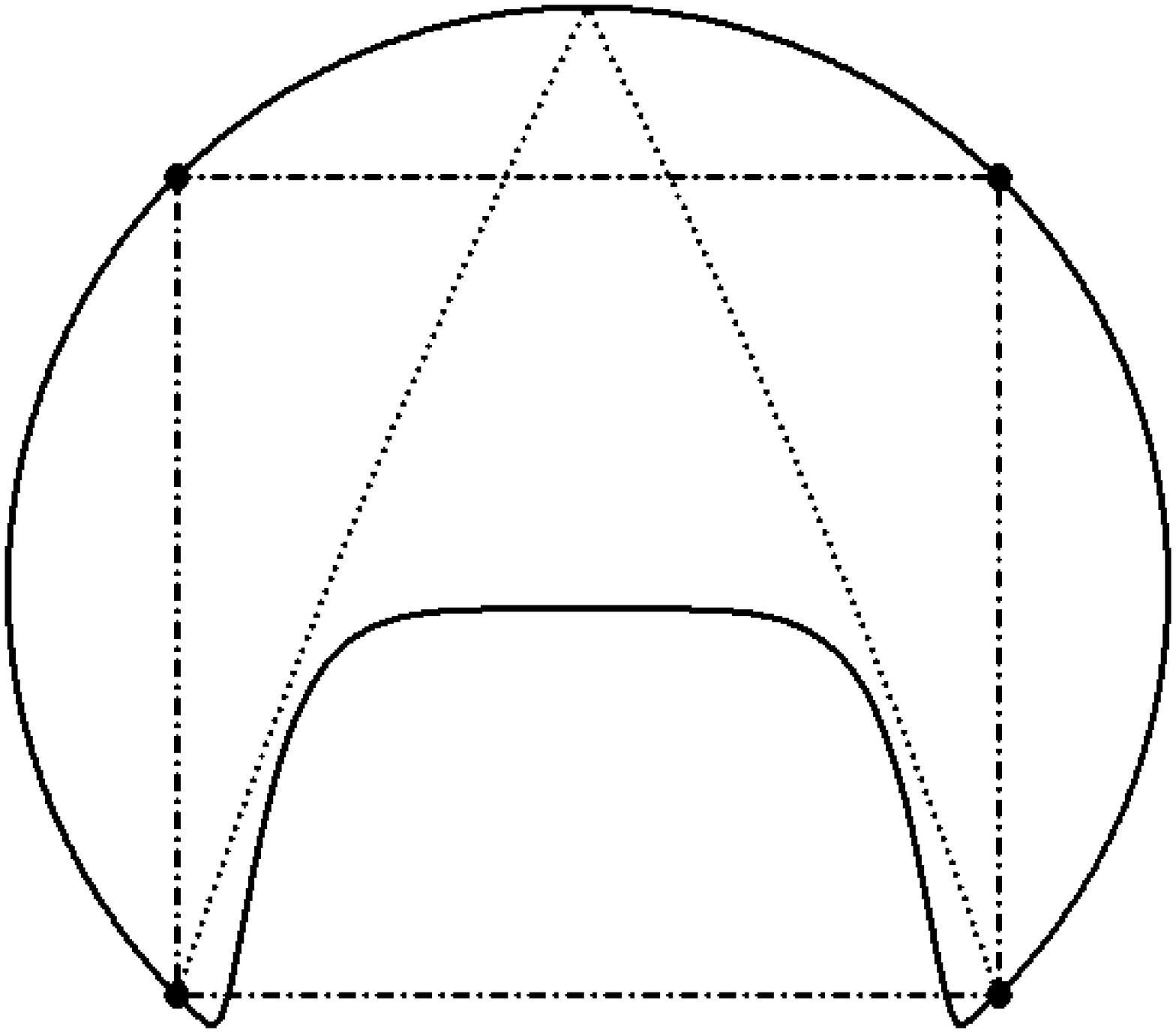}}
\scalebox{.16}{\includegraphics*{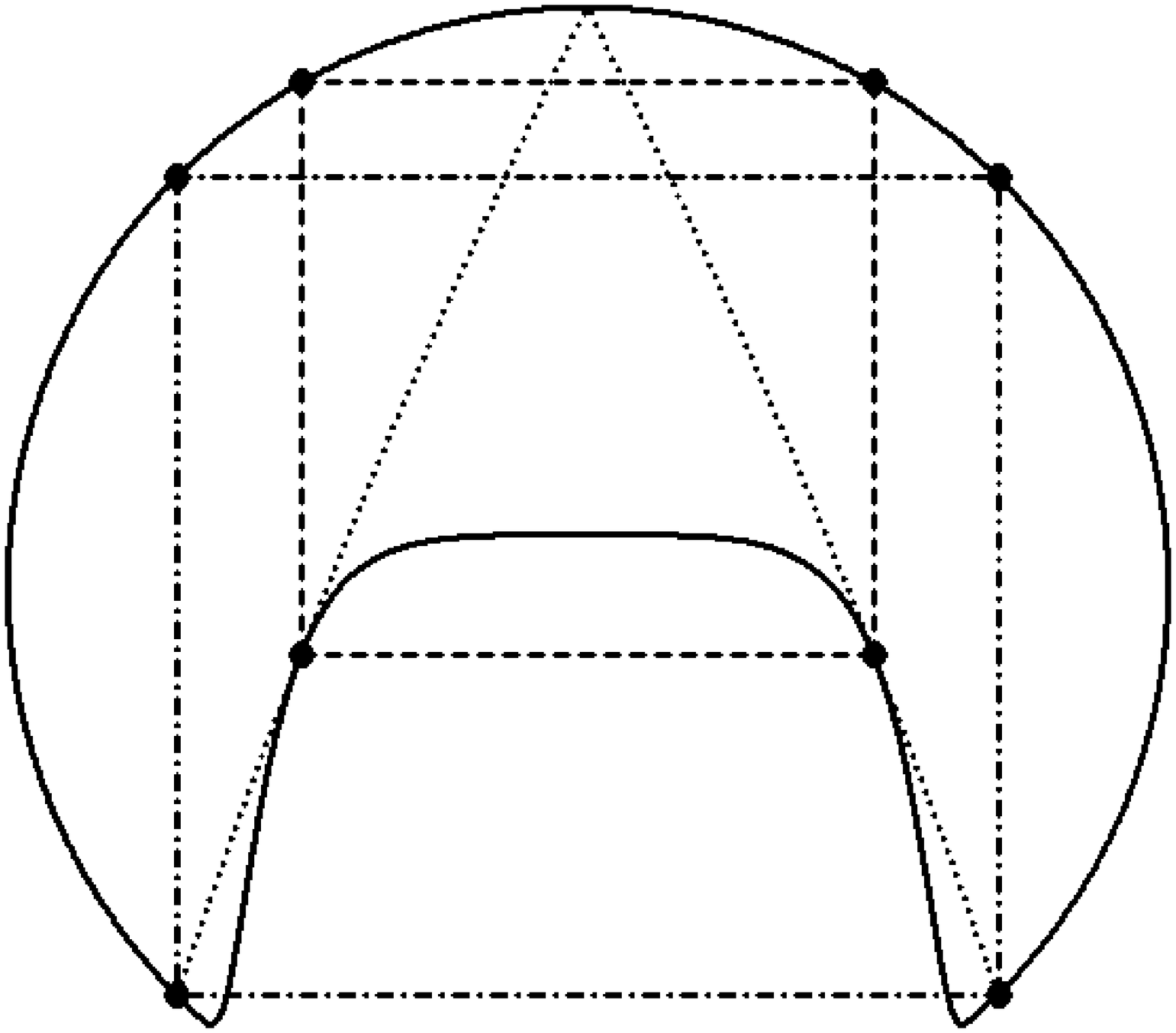}}
  \caption{\small }
  \label{wigglingvsjump}
\end{figure*}

What looks like the equal sides 
(though they are not line segments) 
of an isosceles triangle on that picture 
is the set of points, that are endpoints 
of the base of a square such that the 
top side of the square has endpoints 
that are symmetric about the $y$-axis, 
belong to the unit circle, and have 
$y$-coordinates $\ge {1\over\sqrt2}$. 
(The equation for the two equal sides of 
that triangle is $y={\sqrt{1-x^2}}-2|x|$, ${\frac{-1}{\sqrt{2}}} \le x \le 
{\frac{1}{\sqrt{2}}}$.)

To get a differentiable example we 
replace that arc with the graph of 
$$y= -{\sqrt{1-x^2}} + c\,\exp\Bigl(
-\bigl(
{\frac{.02}{(x+{\frac1{\sqrt{2}}})^2}}+
{\frac{.02}{(x-{\frac1{\sqrt{2}}})^2}}
\bigr)
\Bigr)
\label{Graph} \ \ \ \ \ \ \ (*)
$$ 

Notice that this graph for positive and 
not too big values of $c$ intersects 
the unit circle only at the end-points 
of the arc that was removed. Among these 
values of $c$, for larger $c$ the graph 
intersects each of the two equal sides 
of that isosceles triangle in two points 
(we do not count the endpoints of the 
arc that was removed). For smaller 
values of $c$ the graph does not 
intersect the equal sides of the triangle 
(except at the endpoints of the 
removed arc). Therefore for a certain 
value of $c$ (approximately $1.18264$) on each side of the triangle there is a unique point that belongs to the graph (apart from the endpoint). We sketch 
the proof that the two inscribed 
squares shown on Figure~1, right, 
are the only ones. 

The above considerations show that these 
two squares are the only inscribed squares that have a horizontal side. 
Assume $S$ is an inscribed square 
with no horizontal side. 
If $S$ has three vertices on the 
$\frac34$-circle (i.e. on the union of 
arcs $\w{AB},\w{BC},\w{CD}$, see Fig.2) 
then it follows that the fourth vertex 
would be on the unit circle, on the 
arc that was removed from our curve, 
a contradiction. Let $\w{DA}$ 
denote the graph of $(*)$. Let the 
vertices of $S$ be $E,F,G,H$ (in this 
order) and consider the case when $E,F$ 
belong to the $\frac34$-circle, and 
$G,H$ belong to $\w{DA}$. We only consider two typical cases. 

\begin{figure*} 
 \scalebox{.24}{\includegraphics*{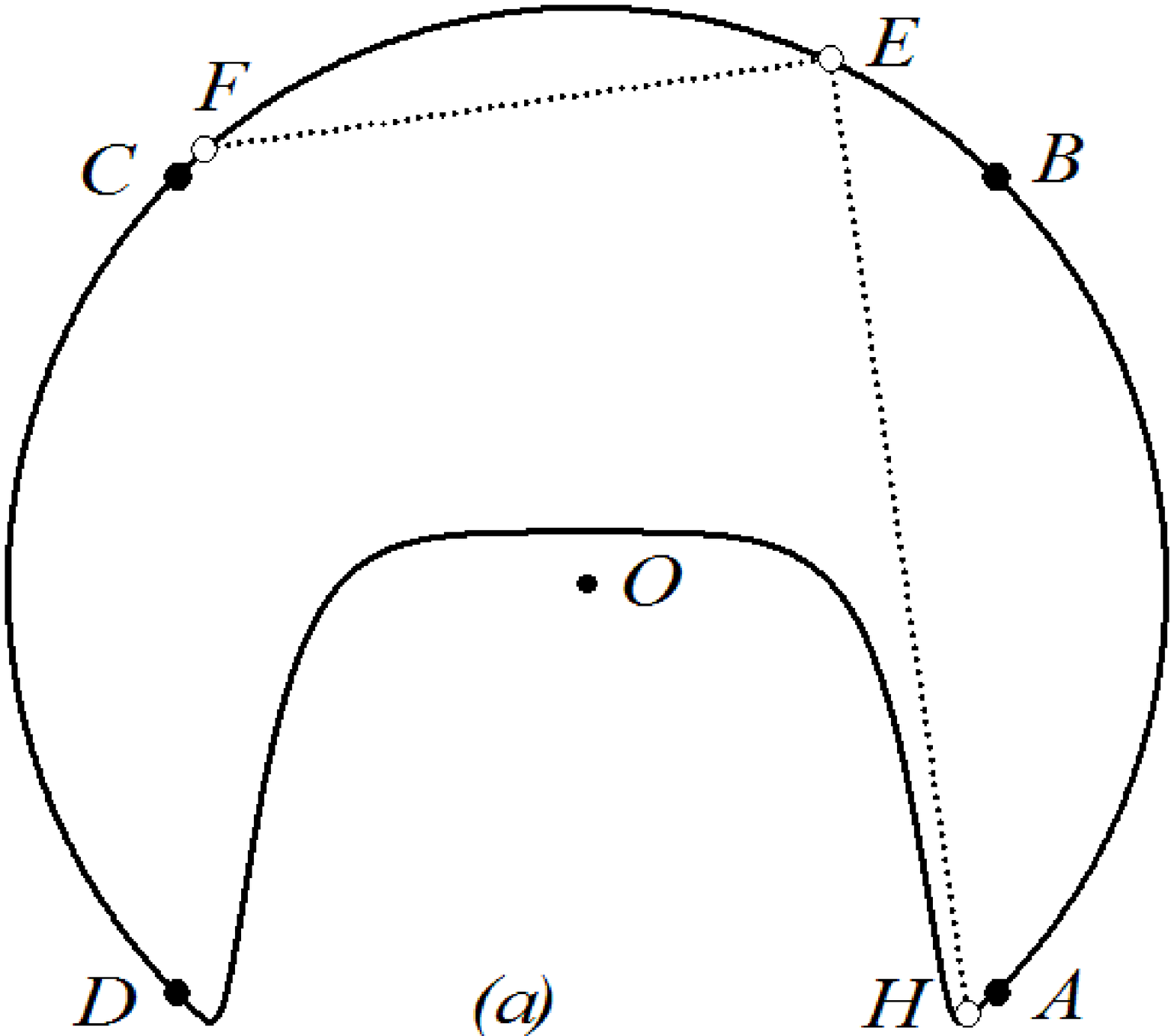}} \hskip12mm \scalebox{.24}{\includegraphics*{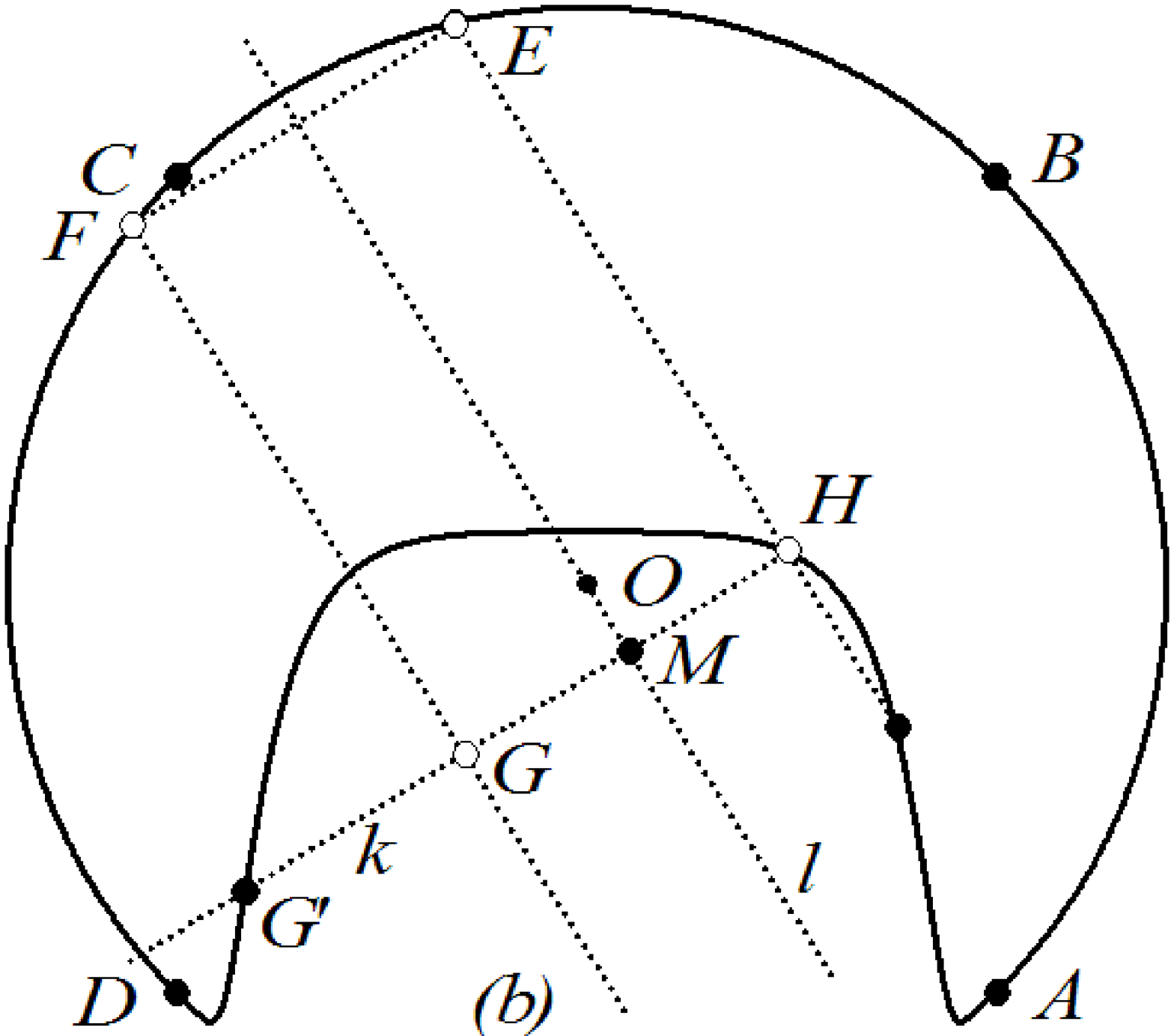}}
  \caption{\small }
  \label{wigglingvsjump}
\end{figure*}

Case 1. $E$ belongs to $\w{BC}$, and 
$H$ belongs to the part of $\w{DA}$ 
that is below the line segment 
$\overline{DA}$ (Fig.2, $(a)$). 
Let $H'$ be the 
intersection of the line through $E,H$ 
with the arc removed from the unit circle. Then $H'$ and $F$ are 
diametrically opposite (since the 
angle at $E$ is right). Hence $F$ 
belongs to $\w{BC}$, and if $d$ denotes 
the distance function, then 
$d(E,F) < \sqrt{2} < d(E,H)$, 
a contradiction. 

Case 2. Assume that $S$ is just 
a rectangle, 
not necessarily a square, and that 
$\overline{FE}$ is a line segment 
with positive slope and endpoints on the 
$\frac34$-circle, and $G,H$ are 
on $\w{DA}$ with $H$ above $\overline{DA}$ (Fig.2, $(b)$). 
Let $l$ be the perpendicular 
bisector of $E$ and $F$. Then $l$ 
goes through the origin $O$ and 
through the midpoint $M$ of $G$ 
and $H$. Let $k$ be the ray starting 
at $M$ and going though $G$, and 
let $G'$ be the ``first'' point on 
$k$ that belongs to $\w{DA}$. 
The reader may verify that 
$d(M,G') > d(M,H)$ and hence 
$G$ does not belong to $\w{DA}$, 
a contradiction. 

\section{How to obtain exactly $n$ squares}  

The next set of examples is also based 
on the idea that we may replace a certain 
arc of the unit circle. It will 
eventually lead to a differentiable 
convex curve with a number of inscribed 
squares specified in advance. 

The idea is very simple, and the proofs 
are easy (though might be technical) 
so we omit some of the details. 

Start with the unit square and this time 
remove the arc 
$[{\frac{-\pi}4},{\frac{\pi}4}]$. 
For convenience we identify any 
real number $P$ with the corresponding 
point on the unit circle, if we treat 
$P$ as an angle. Pick any 
$P \in ({\frac{-\pi}4},{\frac{\pi}4})$ 
and connect $P$ to $A={\frac{-\pi}4}$, 
and $P$ to $B={\frac\pi4}$,
with a circular arc of radius close 
to $1$ but less than $1$. Clearly 
the resulting simple closed curve 
has only two inscribed squares, 
as shown on Fig.3,~$(a)$. 

\begin{figure*} 
 \scalebox{.21}{\includegraphics*{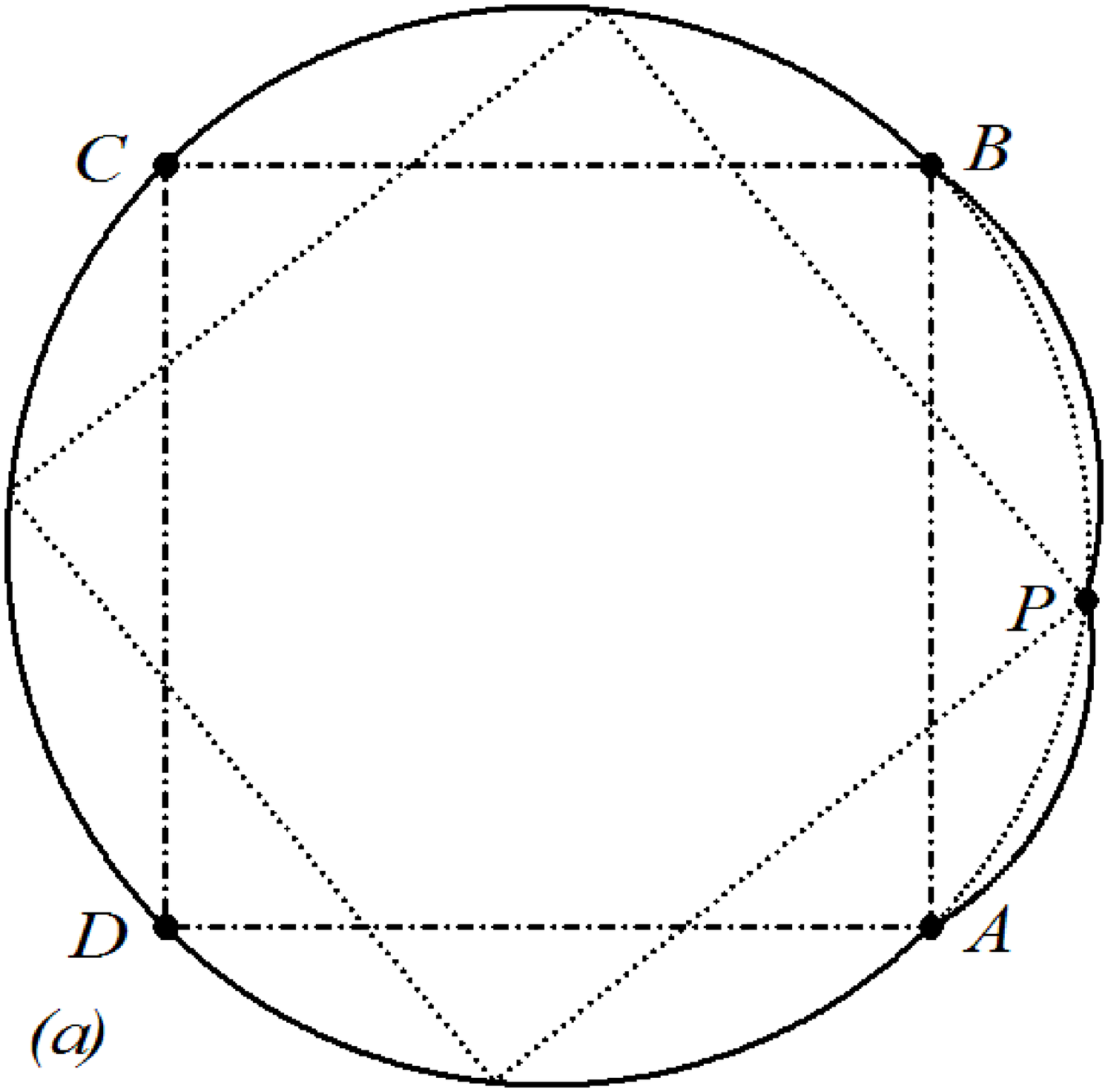}}  \scalebox{.21}{\includegraphics*{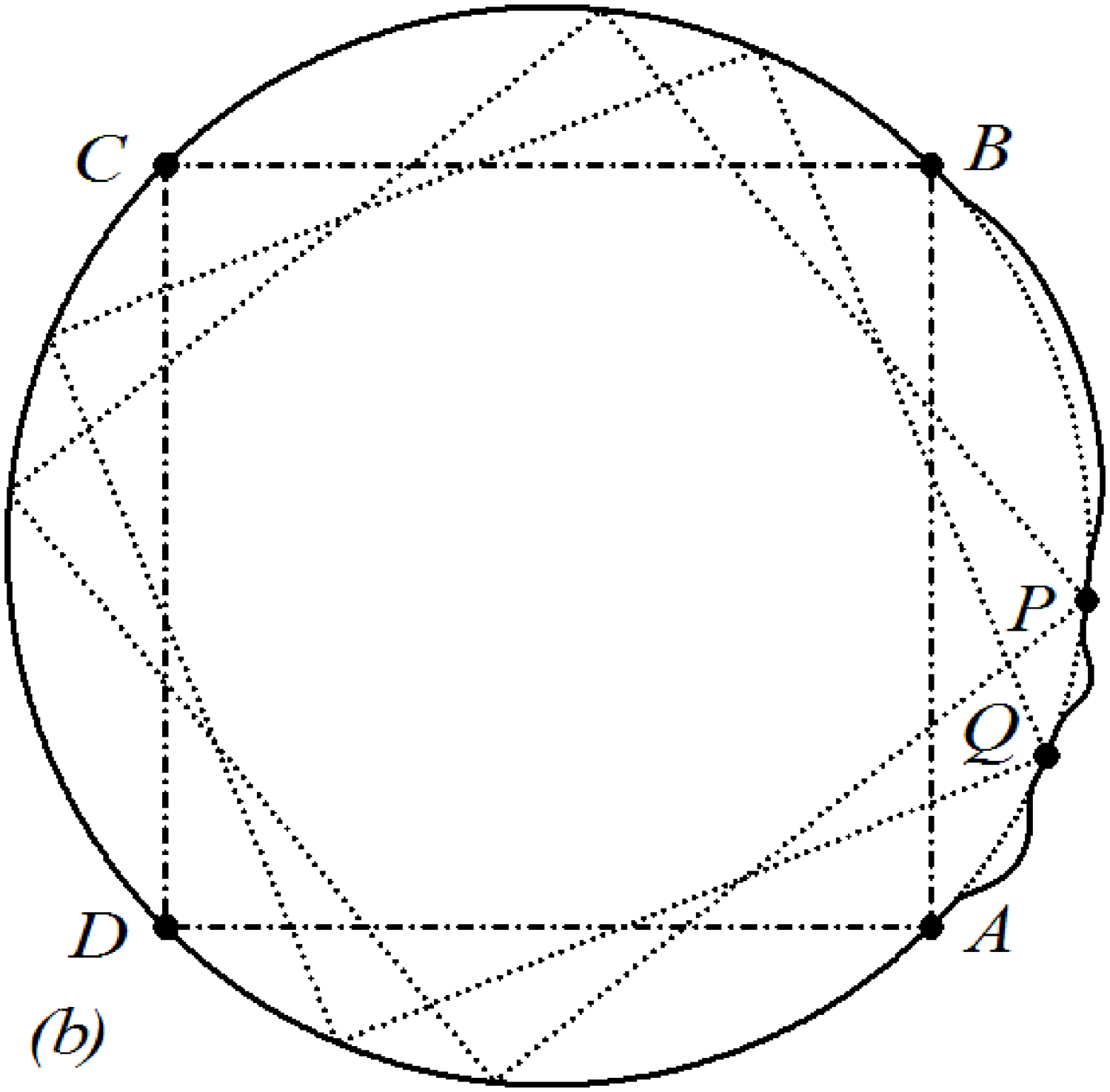}}
\scalebox{.21}{\includegraphics*{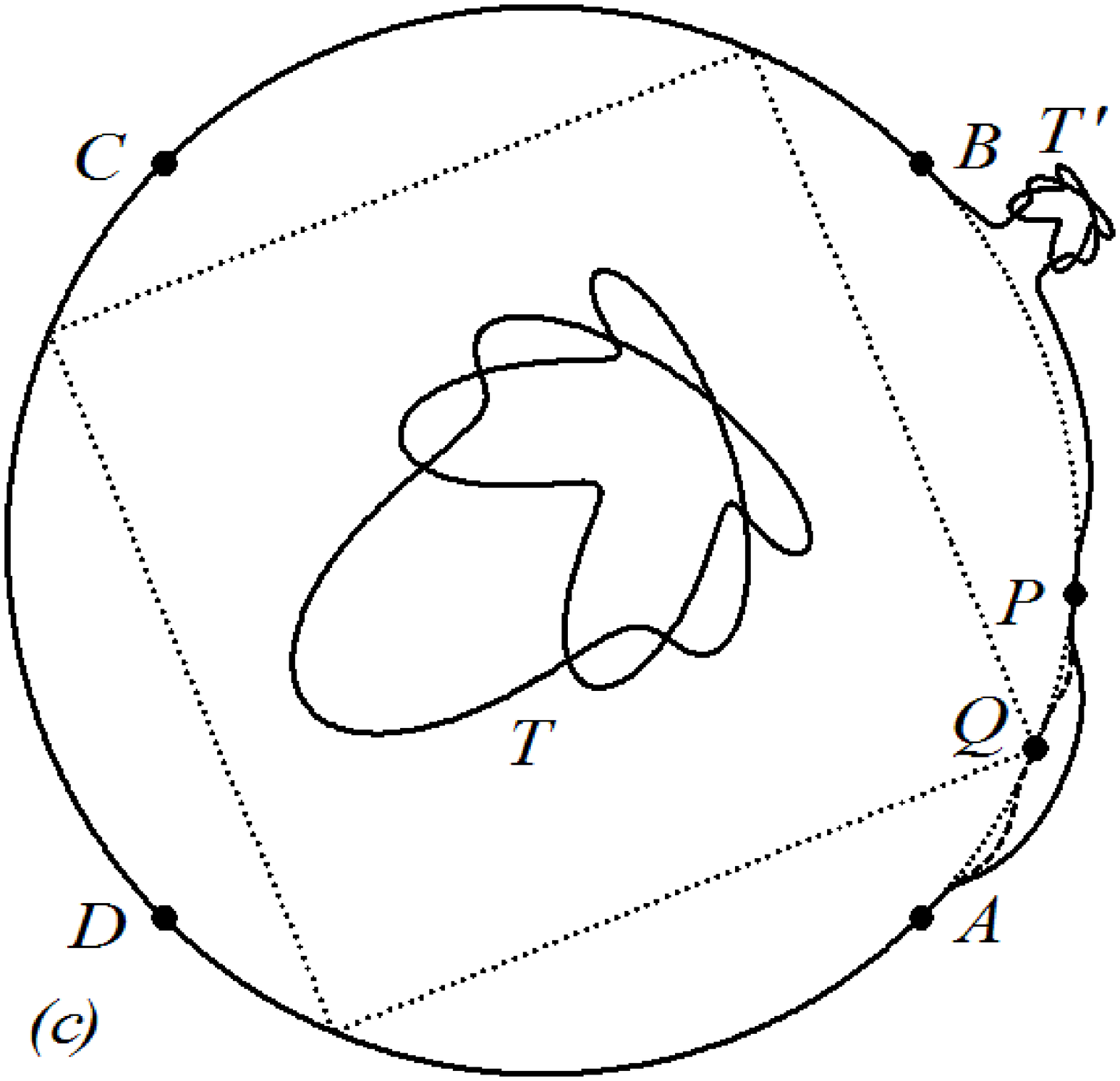}}
  \caption{\small }
  \label{wigglingvsjump}
\end{figure*}

Of course if we add circular arcs of 
smaller than $1$ radius then we do 
not get a differentiable curve, but 
we may instead add an arc of the 
form (polar coordinates): 
$r(\theta) = 
1+c\, \exp
\Bigl(-
\bigl(
{\frac{.02}{(\theta-U)^2}}+
{\frac{.02}{(V-\theta)^2}}
\bigr)
\Bigr)
$ 
in order to connect any given 
pair of points $U$ and $V$ on the 
unit circle, where $c>0$. 
For example, on Fig.3, 
$(b)$, the points $P$ and $B$ are 
connected with an arc of the above 
type, with $c=0.05$. Clearly 
this approach results in an infinitely 
differentiable curve. 
If we select the constant $c$ small 
enough then the signed curvature 
would be positive for all 
$\theta \in [U,V]$, and 
therefore the (region bounded by the) 
simple closed curve obtained in this 
manner would be convex. We can pick 
any finite number of points between 
$A$ and $B$ on the unit circle and replace the consecutive unit circle 
arcs that connect these points with 
arcs of the type described above, and since exactly one inscribed square 
would correspond to each of these 
points we may obtain an infinitely differentiable, convex simple closed 
curve with exactly $n$ inscribed 
squares, for any positive integer $n$ 
given in advance, 
as stated in Theorem~\ref{exactly-n}. 
See Fig.2, $(b)$. 

\section{On the role of $St$, $J^+$ and $J^-$} 

In this section we indicate a 
possible proof of the following conjecture. 

\begin{conjecture} 
Given any immersed curve $T$ in 
the plane, there is a positive 
integer $m$ such that for every 
$n \ge m$ there is an immersed curve 
$T_n$ which has the same values of 
$St$, $J^+$ and $J^-$ as $T$, and such 
that $T_n$ has exactly $n$ inscribed 
squares. Moreover there is $k$ 
(independent of $n$) such that 
all but $k$ many of the inscribed 
squares of $T_n$ have the property 
that their vertices appear in the 
same order in which they appear on 
$T_n$. 
\end{conjecture}

The idea is the following. Start 
with an immersed curve $T$ (e.g. 
the one shown in Fig.3, $(c)$, 
in the middle of the circle). Pull 
one of the loops of $T$ and wrap 
it around the unit circle, and at the 
same time make the rest of $T$ much smaller, so that we have a very small 
copy of $T$, very close to $B$, 
as shown in Fig.3, $(c)$, except for 
the loop that is wrapped around the 
unit circle. Call the resulting curve 
$T'$. More precisely we assume that 
a point 
$P \in ({\frac{-\pi}4},{\frac{\pi}4})$ 
has been fixed, the unit circle 
arcs from $A$ to $P$, and from $P$ to $B$ have been replaced by arcs of the type described above, and then $T'$ has been formed by wrapping one of its loops around, so that, except for this loop, a very small (topological) copy 
of $T$ remains very close to $B$, and 
``between'' $P$ and $B$. 
We also assume that $T'$ 
is differentiable. 

Clearly $T$ and $T'$ have the same 
values for $St$, $J^+$ and $J^-$. 
We will give a proof of our conjecture 
based on the following genericity 
assumption (GA), which we leave without 
proof. (We do not know how to prove it, but we believe it is correct.) 

\begin{ga} 
The above transformation of $T$ 
to $T'$ can be done in such a way 
that $T'$ has only finitely many 
inscribed squares. 
\end{ga} 

Now in order to prove our conjecture 
based on our GA, let $m$ be the finite 
number of inscribed squares of $T'$. 
Let $k$ be the number of them for 
which the vertices  
appear in order different from 
the order in which they appear on $T'$. 
We can pick points $Q$ on the unit 
circle, between $A$ and $P$, one at 
a time, replacing an arc of the type 
described above with two smaller arcs, 
so that every time a new inscribed square 
with one vertex at the new point $Q$ 
would be introduced, and no other 
inscribed squares would be introduced. 
Notice that the new squares have vertices 
which appear in the same order as 
in $T'$ (see Fig.3, $(c)$). 
This completes the proof.

\bibliographystyle{amsplain}

\end{document}